\documentclass[a4paper]{amsart}
\usepackage[utf8]{inputenc}
\usepackage[T1]{fontenc}
\pdfoutput=1

\usepackage[style=alphabetic,maxbibnames=10,isbn=false,backend=bibtex]{biblatex}
\DeclareSourcemap{
	\maps[datatype=bibtex, overwrite]{
		\map{
			\step[fieldset=pagetotal, null]
			\step[fieldset=address, null]
		}
	}
}
\AtEveryBibitem{%
  \clearfield{pagetotal} 
  \clearlist{address}
  \clearlist{location} 
}
\appto{\bibsetup}{\emergencystretch=1em}
\usepackage{paralist}
\usepackage{hyperref}

\makeatletter
\newcommand{\proofstep}[1]{%
  \par%
  \addvspace{\medskipamount}%
  \textit{#1\@addpunct{.}}\enspace\ignorespaces
}  
  
\newcommand{\doublerightarrow}{%
  \mathrel{%
  \mathop{\vcenter{%
    \offinterlineskip\ialign{##\cr
    $\rightarrow$\cr\noalign{\kern.5ex}
    $\rightarrow$\cr}}}}
}
\newcommand*{\xdoublerightarrow}[2]{\mathrel{%
  \settowidth{\@tempdima}{$\scriptstyle#1$}
  \settowidth{\@tempdimb}{$\scriptstyle#2$}
  \ifdim\@tempdimb>\@tempdima \@tempdima=\@tempdimb\fi
  \mathop{\vcenter{%
    \offinterlineskip\ialign{\hbox to\dimexpr\@tempdima+1em{##}\cr
    \rightarrowfill\cr\noalign{\kern.5ex}
    \rightarrowfill\cr}}}\limits^{\!#1}_{\!#2}}}%
    
\newcommand*{\xmanyrightarrows}[2]{\mathrel{%
  \settowidth{\@tempdima}{$\scriptstyle#1$}
  \settowidth{\@tempdimb}{$\scriptstyle#2$}
  \ifdim\@tempdimb>\@tempdima \@tempdima=\@tempdimb\fi
  \mathop{\vcenter{%
    \offinterlineskip\ialign{\hbox to\dimexpr\@tempdima+1em{##}\cr
    \rightarrowfill\cr
    \hfill\vdots\hfill\cr\noalign{\kern.5ex}
    \rightarrowfill\cr}}}\limits^{\!#1}_{\!#2}}}%
    
\makeatother
	
\newenvironment{acknowledgement}%
  {\par%
   \hspace{\parindent}\textsc{Acknowledgments.}\enskip\ignorespaces}
  {\par}

\usepackage{amsmath, amssymb, amsfonts, mathtools, amsthm, mathrsfs, dsfont, fge, extpfeil}
\usepackage{tikz,tikz-cd,multirow,faktor,xfrac}
\usetikzlibrary{matrix,arrows,positioning,calc}

\newcommand\setC{\mathbb{C}}
\newcommand\F{\mathbb{F}}

\newcommand\Q{\mathbb{Q}}
\newcommand\Z{\mathbb{Z}}
\newcommand{\setN}{\mathbb{N}}
\newcommand{\setP}{\mathbb{P}}

\newcommand\setH{\mathbb{H}}

\newcommand\B{\mathcal{B}}
\newcommand\C{\mathcal{C}}
\newcommand\M{\mathcal{M}}
\newcommand\N{\mathcal{N}}

\newcommand\cR{\mathcal{R}}
\newcommand\cP{\mathcal{P}}
\newcommand\cQ{\mathcal{Q}}

\newcommand\del{\partial}
\newcommand\eps{\varepsilon}

\newcommand\id{\mathrm{id}}
\DeclareMathOperator{\rad}{rad}

\DeclareMathOperator{\im}{im}
\DeclareMathOperator{\coker}{coker}
\DeclareMathOperator{\Ker}{Ker}
\def\embeds{\hookrightarrow}
\DeclareMathOperator{\onto}{\twoheadrightarrow}
\def\setminus{\mathbin{\fgebackslash}}
\DeclareMathOperator{\tr}{tr}

\DeclareMathOperator{\Char}{char}

\DeclareMathOperator{\mods}{mod}

\DeclareMathOperator{\add}{add}
\DeclareMathOperator{\SL}{SL}
\DeclareMathOperator{\Sl}{\mathcal{S}\ell}
\DeclareMathOperator{\GL}{GL}
\DeclareMathOperator{\Hom}{\mathrm{Hom}}
\DeclareMathOperator{\End}{\mathrm{End}}
\DeclarePairedDelimiter{\Euler}{\langle}{\rangle}

\DeclarePairedDelimiter{\card}{|}{|}
\DeclarePairedDelimiter{\norm}{||}{||}

\DeclarePairedDelimiter{\ceil}{\lceil}{\rceil}

\def\Ext#1#2{\mathrm{Ext}_{#1}^{#2}}
\def\defined#1{\textbf{#1}}

\def\isom{\mathrel{\cong}}

\def\spann#1{\left\langle #1 \right\rangle}

\DeclareMathOperator{\adj}{adj}

\theoremstyle{definition}
\newtheorem{mydef}{Definition} 

\theoremstyle{plain}
\newtheorem{theorem}[mydef]{Theorem}
\newtheorem*{theorem*}{Theorem}

\newtheorem{prop}[mydef]{Proposition}
\newtheorem{proposition}[mydef]{Proposition}
\newtheorem{lemma}[mydef]{Lemma}
\newtheorem{corollary}[mydef]{Corollary}
\newtheorem*{corollary*}{Corollary}

\theoremstyle{remark}

\newtheorem*{example*}{Example}

\newtheorem*{remark}{Remark}

\newtheorem*{remarks}{Remarks}

\title{Wild Kronecker quivers and amenability}
\author{Sebastian Eckert}
\address{Sebastian Eckert, Fakultät für Mathematik, Universität Bielefeld, Postfach 100~131, 33501~Bielefeld, Germany.} \email{seckert@math.uni-bielefeld.de}
\curraddr{Max-Planck-Institut für Mathematik, Vivatsgasse 7, 53111~Bonn, Germany.} \email{eckert@mpim-bonn.mpg.de}
\thanks{During his PhD studies, the author was supported~by the Alexander von Humboldt Foundation in the framework of an Alexander von Humboldt Professorship endowed by the German Federal Ministry of Education and Research.
An updated version of these notes was later prepared during a stay at the Max Planck Institute for Mathematics in Bonn, and the authors is grateful for its hospitality and financial support.
}


\begin{document}

\begin{abstract}
We apply the notion of hyperfinite families of modules to the wild path algebras of generalised Kronecker quivers $k\Theta(d)$. While the preprojective and postinjective component are hyperfinite, we show the existence of a family of non-hyperfinite modules in the regular component for some $d$. Making use of dimension expanders to achieve this, our construction is more explicit than previous results.
From this it follows that no finitely controlled wild algebra is of amenable representation type.
\end{abstract}

\maketitle

\section{Introduction}
The notions of hyperfiniteness for countable sets of modules and amenable representation type for algebras have been introduced by \cite{Elek2017InfiniteDimensionalRepresentationsAmenabilty}.
We will work with the definitions as follow.

\begin{mydef} \label{def:HyperfinitenessAmenability}
Let~$k$ be a field, $A$ be a finite dimensional $k$-algebra and let $\M$ be a set of finite dimensional $A$-modules. One says that $\M$ is \defined{hyperfinite} provided for every $\varepsilon > 0$ there exists a number $L_\varepsilon > 0$ such that for every $M \in \M$ there exist both, a submodule $N \subseteq M$ such that
\begin{equation} \label{eq:HFSubmoduleBig} \dim_k N \geq (1-\varepsilon) \dim_k M, \end{equation} and modules $N_1,N_2, \dots N_t \in \mods A$, with $\dim_k N_i \leq L_\varepsilon$, such that $N \isom \bigoplus_{i=1}^{t} N_i$.

The $k$-algebra $A$ is said to be of \defined{amenable representation type} provided the set of all finite dimensional $A$-modules (or more precisely, a set which meets every isomorphism class of finite dimensional $A$-modules) is hyperfinite.
\end{mydef}

Previously, the author has shown that tame quiver algebras are of amenable representation type (giving a new proof and not using a previous result on string algebras) while wild quiver algebras are not (using results of Elek), thus working towards \cite[Conjecture~1]{Elek2017InfiniteDimensionalRepresentationsAmenabilty}, presuming that finite dimensional algebras are of tame type if and only if they are of amenable representation type. 

In this note we will focus on hyperfinite families for (wild) Kronecker algebras while also ascertaining that wild Kronecker algebras and more generally, finitely controlled wild algebras are not of amenable representation type.

We further use the following facts from \cite{Eckert2019}.

\begin{proposition} \label{prop:AdditiveClosureStaysHyperfinite}
Let $\M$ be a family of $A$-modules. If $\M$ is hyperfinite, so is the family of all finite direct sums of modules in $\M$.
\end{proposition}

\begin{proposition} \label{prop:ExtendingHFfromSubmodulesOfBoundedCodimension} 
Let $A$ be a finite dimensional $k$-algebra. Let $\M,\N \subset \mods A$ where $\N$ is hyperfinite. If there is some $L \geq 0$ such that for all $M \in \M$, there exists a submodule $P \subset M$ of codimension less than or equal to $L$, and $P \in \N$, then $\M$ is also hyperfinite.
\end{proposition}

\begin{proposition} \label{prop:HFPreservingFunctors} 
Let $k$ be a field and $A,B$ be two finite dimensional $k$-algebras. Let $F \colon \mods A \to \mods B$ be an additive, left-exact functor such that there exists $K_1, K_2 > 0$ with
\begin{equation}  \label{eq:EquivalenceCondHFPresFun}
K_1 \dim X \leq \dim F(X) \leq K_2 \dim X,\end{equation}
 for all $X \in \mods A$.
If $\N \subseteq \mods A$ is a hyperfinite family, then the family $\{F(X) \colon X \in \N\} \subseteq \mods B$ is also hyperfinite.
\end{proposition}

\section{The special case of the $2$-Kronecker quiver} \label{section:KroneckerQuiver}

Let us first recall the situation for the tame 2-Kronecker quiver. It follows from the results on string algebras in \cite[Proposition~10.1]{Elek2017InfiniteDimensionalRepresentationsAmenabilty}, but we give a direct and independent proof here for illustration purposes and to the convenience of the reader.

\begin{lemma} \label{lemma:DescentOnDefectOfPIs}
Let $X = \tau^{r} I(i)$ be some indecomposable postinjective $kQ$-module of defect $\del(X) = d$. Then there is an injective module $I(j)$ such that there exists a non-zero morphism $\theta \colon X \to I(j)$. Moreover,
for any direct summand $Z$ of $\ker \theta$, we have $\del(Z) < d$.
\end{lemma}

\begin{theorem} \label{thm:2KroneckerAmenable}
Let $k$ be any field. Then the path algebra of the $2$-Kronecker quiver $\Theta(2)$ is of amenable representation type.
\begin{proof}
We fix notation for the vertices and arrows as $1 \xdoublerightarrow{a}{b} 2$.

It is well-known (see, e.g., \cite[Theorem~4.3.2]{Benson1998RepresentationsCohomologyI} or \cite{Burgermeister1986ClassificationRepresentationsDoubleFleche}) that the indecomposable preprojective and postinjective $k$-representations of $Q$ are parametrised by 
\[ P_n \colon k^{n} \xdoublerightarrow{\begin{bsmallmatrix} \id \\ 0\end{bsmallmatrix}}{\begin{bsmallmatrix} 0 \\ \id\end{bsmallmatrix}} k^{n+1}, \quad \text{and} \quad
Q_n \colon k^{n+1} \xdoublerightarrow{\begin{bsmallmatrix} \id & 0\end{bsmallmatrix}}{\begin{bsmallmatrix} 0 & \id\end{bsmallmatrix}} k^{n},\;\text{respectively},\]

both for $n\geq 0$, while the indecomposable regular representations can be parametrised by
\[k^n \xdoublerightarrow{\phi}{\psi} k^n,\]

where either $\phi$ is the identity and $\psi$ is given by the companion matrix of a power of a monic irreducible polynomial over $k$, or $\psi$ is the identity and $\phi$ is given by the companion matrix of a monomial. 

We will show that the preprojective component $\cP$, the regular component $\cR$ and the postinjective component $\cQ$ are each hyperfinite families to conclude the amenability of $\mods kQ$. We will give an argument for the indecomposable objects in each component and then apply Proposition~\ref{prop:AdditiveClosureStaysHyperfinite} to extend the result.

\begin{figure}[ht]
\centering
\begin{tikzpicture}[mystyle/.style={insert path={\pgfextra{\node [right] at (\tikzlastnode.north) {\tiny\tikzlastnode};}}},level distance=3em,sibling distance=3em,edge from parent = [->,-stealth]]

\def\n{10}
\def\m{5}
\def\restart{8}
\def\exampleK{3}
\pgfmathparse{int(\restart+1)} \let\restartt\pgfmathresult

\node (e-1) {$e_1$}
    child {node (f-1) {$f_1$} edge from parent[blue]}
    child {node (f-2) {$f_2$} edge from parent[orange]};

\foreach [evaluate=\number as \prev using int(\number-1),evaluate=\number as \succ using int(\number+1),evaluate={\thirdle=int(Mod(\number,\exampleK)}] \number in {2,...,\m}{%
	\if \thirdle0
		\node[draw=red,thick,right of=e-\prev] (e-\number) {$e_\number$}
		child[missing]
		child {node (f-\succ) {$f_\succ$} edge from parent[orange]};
	\else
		\node[right of=e-\prev] (e-\number) {$e_\number$}
		child[missing]
		child {node (f-\succ) {$f_\succ$} edge from parent[orange]};
	\fi
}

\node[right of=e-\m] (e-x) {$\dots$}
    child[missing]
    child {node (f-x) {$\dots$} edge from parent[draw=none]};
\node[right of=e-x] (e-y) {$\dots$}
    child[missing]
    child[missing];
    
\pgfmathparse{int(\n-\restart)} \let\mm\pgfmathresult    
\pgfmathparse{int(Mod(\restart,\exampleK))} \let\thirdle\pgfmathresult
\if \thirdle0
	\node[draw=red,thick,right of=e-y] (e-\restart) {$e_{n-\m}$}
\else
	\node[right of=e-y] (e-\restart) {$e_{n-\mm}$}
\fi
    child {node (f-\restart) {$f_{n-\mm}$} edge from parent[blue]}
    child[missing];

\pgfmathparse{int(\n-1)} \let\nn\pgfmathresult
\pgfmathparse{int(\n+1)} \let\nnn\pgfmathresult
\foreach [evaluate=\number as \prev using int(\number-1),evaluate={\thirdle=int(Mod(\number,\exampleK)}] \number in {\restartt,...,\nn}{%
	\pgfmathparse{int(\n-\number)} \let\index\pgfmathresult
	\if \thirdle0
		\node[draw=red,thick,right of=e-\prev] (e-\number) {$e_{n-\index}$}
		child {node (f-\number) {$f_{n-\index}$} edge from parent[blue]}
		child[missing]
	\else
		\node[right of=e-\prev] (e-\number) {$e_{n-\index}$}
		child {node (f-\number) {$f_{n-\index}$} edge from parent[blue]}
	\fi
		child[missing];
}    

\pgfmathparse{int(Mod(\n,\exampleK))} \let\thirdle\pgfmathresult
\if \thirdle0
	\node[draw=red,thick,right of=e-\nn] (e-\n) {$e_{n}$}
\else
	\node[right of=e-\nn] (e-\n) {$e_{n}$}
\fi
		child {node (f-\n) {$f_{n}$} edge from parent[blue]}
		child {node (f-\nnn) {$f_{n+1}$} edge from parent[orange]};

\path[->,-stealth]
	\foreach [evaluate={\thirdle=int(Mod(\number,\exampleK)}] \number in {2,...,\m}{
		\if \thirdle0
			(e-\number) edge[blue] (f-\number)
		\else
			(e-\number) edge[blue] (f-\number) 
		\fi
	}
	\foreach [evaluate=\number as \succ using int(\number+1),evaluate={\thirdle=int(Mod(\number,\exampleK)}] \number in {\restart,...,\nn}{
		\if \thirdle0
			(e-\number) edge[orange] (f-\succ)
		\else
			(e-\number) edge[orange] (f-\succ)
		\fi
	};
\end{tikzpicture}
\caption{\label{fig:VisualiseProof2KroneckerPP}The coefficient quiver of $P_n$, showing a decomposition for $n \equiv 1 \mod 3$ and $K_\eps = 3$.}
\end{figure}
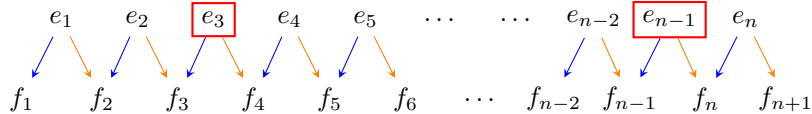

We start with the preprojectives and let $\varepsilon > 0$. Set $K_\varepsilon \coloneqq \ceil{\frac{1}{2\varepsilon}} + 1$ and $L_\varepsilon = \frac{1}{\varepsilon} + 3$.
Let $X = P_n$ be some indecomposable preprojective. If $\dim X \leq L_\varepsilon$, there is nothing to show. We may thus assume that $\dim X > L_\varepsilon$, implying $n \geq K_\varepsilon$, and write $n = j \cdot K_\varepsilon + r,$ where $0 \leq r < K_\varepsilon$. Now consider the standard basis $\{e_1, e_2, \dots e_n\}$ of $k^n$. Let $U$ be the submodule of $X$ generated by the subset
\begin{align*} & \{e_1, \dots, e_{K_\varepsilon-1}\} \cup \{e_{K_\varepsilon+1}, \dots, e_{2 K_\varepsilon - 1}\} \cup \dots \\ 
\cup & \{e_{(j-1) K_\varepsilon + 1}, \dots, e_{j K_\varepsilon - 1}\} \cup \{ e_{j K_\varepsilon + 1}, \dots, e_{n} \},\end{align*} 
dropping  every $K_\varepsilon$-th basis vector. Then $U$ decomposes into $j$ direct summands of type $P_{K_\varepsilon-1}$ and a smaller rest term in case $r\neq 0$. All summands will thus be of $k$-dimension smaller than $2 (K_\varepsilon - 1) + 1 < L_\varepsilon$. Moreover,
\begin{align*}
\dim U &= \dim X - j  = \dim X - \frac{n - r}{K_\varepsilon} = \dim X - \frac{\dim X -1}{2 K_\varepsilon} + \frac{r}{K_\varepsilon} \\
& \geq \dim X - \varepsilon (\dim X - 1) > (1-\varepsilon) \dim X.
\end{align*}
This shows that the family of indecomposable preprojective modules $\cP(kQ)$ is hyperfinite. We exemplify this process in Figure~\ref{fig:VisualiseProof2KroneckerPP}.

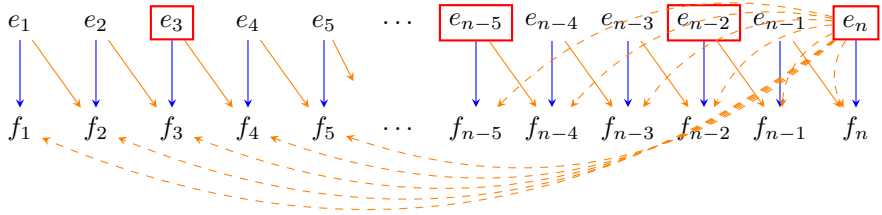
\begin{figure}[ht]
\centering
\begin{tikzpicture}[mystyle/.style={insert path={\pgfextra{\node [right] at (\tikzlastnode.north) {\tiny\tikzlastnode};}}},level distance=4em,sibling distance=3em,edge from parent = [->,-stealth]]
\def\n{14}
\def\m{5}
\def\restart{9}
\def\exampleK{3}
\pgfmathparse{int(\restart+1)} \let\restartt\pgfmathresult

\node (e-1) {$e_1$}
    child {node (f-1) {$f_1$} edge from parent[blue]};

\foreach [evaluate=\number as \prev using int(\number-1),evaluate=\number as \succ using int(\number+1),evaluate={\thirdle=int(Mod(\number,\exampleK)}] \number in {2,...,\m}{%
	\if \thirdle0
		\node[draw=red,thick,right of=e-\prev] (e-\number) {$e_\number$}
		child {node (f-\number) {$f_\number$} edge from parent[blue]};
	\else
		\node[right of=e-\prev] (e-\number) {$e_\number$}
		child {node (f-\number) {$f_\number$} edge from parent[blue]};
	\fi
}

\node[right of=e-\m] (e-x) {$\dots$}
    child {node (f-x) {$\dots$} edge from parent[draw=none]};
    
\pgfmathparse{int(\n-\restart)} \let\mm\pgfmathresult    
\pgfmathparse{int(Mod(\restart,\exampleK))} \let\thirdle\pgfmathresult
\if \thirdle0
	\node[draw=red,thick,right of=e-x] (e-\restart) {$e_{n-\m}$}
	    child {node (f-\restart) {$f_{n-\mm}$} edge from parent[blue]};
\else
	\node[right of=e-x] (e-\restart) {$e_{n-\mm}$}
	    child {node (f-\restart) {$f_{n-\mm}$} edge from parent[blue]};
\fi

\pgfmathparse{int(\n-1)} \let\nn\pgfmathresult
\pgfmathparse{int(\n+1)} \let\nnn\pgfmathresult
\foreach [evaluate=\number as \prev using int(\number-1),evaluate={\thirdle=int(Mod(\number,\exampleK)}] \number in {\restartt,...,\nn}{%
	\pgfmathparse{int(\n-\number)} \let\index\pgfmathresult
	\if \thirdle0
		\node[draw=red,thick,right of=e-\prev] (e-\number) {$e_{n-\index}$}
		child {node (f-\number) {$f_{n-\index}$} edge from parent[blue]};
	\else
		\node[right of=e-\prev] (e-\number) {$e_{n-\index}$}
		child {node (f-\number) {$f_{n-\index}$} edge from parent[blue]};
	\fi
}    

\pgfmathparse{int(Mod(\n,\exampleK))} \let\thirdle\pgfmathresult
	\node[draw=red,thick,right of=e-\nn] (e-\n) {$e_{n}$}
		child {node (f-\n) {$f_{n}$} edge from parent[blue]};

\path[->,-stealth]
	\foreach [evaluate=\number as \succ using int(\number+1),evaluate={\thirdle=int(Mod(\number,\exampleK)}] \number in {1,...,\m}{
		\if \thirdle0
			(e-\number) edge[orange] (f-\succ)
		\else
			(e-\number) edge[orange] (f-\succ) 
		\fi
	}
	\foreach [evaluate=\number as \succ using int(\number+1),evaluate={\thirdle=int(Mod(\number,\exampleK)}] \number in {\restart,...,\nn}{
		\if \thirdle0
			(e-\number) edge[orange] (f-\succ)
		\else
			(e-\number) edge[orange] (f-\succ)
		\fi
	};
\path[->,-stealth]
	\foreach \number in {1,...,\m}{(e-\n) edge[orange,dashed,bend left] (f-\number)}
	\foreach \number in {\restart,...,\n}{(e-\n) edge[orange,dashed,bend right] (f-\number)};
\end{tikzpicture}
\caption{\label{fig:VisualiseProof2KroneckerRegular}The coefficient quiver of some $R_n(\id,\varphi)$, exhibiting a way to find a suitable submodule for $n \equiv 2 \mod 3$, $K_\eps = 3$.}
\end{figure}

If $X = R_n(\phi,\psi)$ is an indecomposable regular module, we may consider the submodule~$Y$ generated by the basis vectors $\{e_1, \dots, e_{n-1}\}$ of the vector space at vertex~$1$. Note that we assume that $\psi$ corresponds to the Frobenius companion matrix of a power of a monic polynomial. Then $Y \isom P_{n-1}$, so by the above it belongs to the hyperfinite family $\cP$.
We have that $\dim Y  = \dim X - 1$. Thus, an application of Proposition~\ref{prop:ExtendingHFfromSubmodulesOfBoundedCodimension} with $H=1$ gives the hyperfiniteness of the indecomposable regular modules.
See Figure~\ref{fig:VisualiseProof2KroneckerRegular} for an example.

We are left to deal with the postinjective case. By Lemma~\ref{lemma:DescentOnDefectOfPIs}, for each indecomposable postinjective module $X$, we can find a submodule $Y \coloneqq \ker \theta$ of strictly smaller defect. Moreover, if $Y$ had a postinjective summand $Z$, it must have defect $\del(Z) < \del(X)$.
In this situation, all indecomposable postinjective modules have defect $d=1$. Choose the hyperfinite family $\N_0 = \cP \cup \cR$ of all preprojective and regular modules. For all postinjective indecomposables, the submodule $Y$ must be in $\add \N_0$, since there are no non-zero postinjective modules $Z$ with defect $\del(Z) < 1$. This family is hyperfinite by the above.
Moreover, the codimension of $Y$ is bounded by the dimension of the indecomposable injectives, of which there are only two. Hence, we can use Proposition~\ref{prop:ExtendingHFfromSubmodulesOfBoundedCodimension} to prove the hyperfiniteness of the indecomposable postinjectives.

Now apply Proposition~\ref{prop:AdditiveClosureStaysHyperfinite} to $\cP \cup \cR \cup \cQ$ to see that $\mods kQ$ is hyperfinite, and thus $kQ$ is of amenable representation type.
\end{proof}
\end{theorem}

\section{Hyperfiniteness from Fragmentability and for exceptional modules}

We will continue by considering the path algebra $k\Theta(d)$ of a generalised, wild ($d\geq3$) Kronecker quiver and show that the indecomposable preprojective and postinjective modules for these algebras form hyperfinite families.
We start with a result connecting to the notion of fragmentability from graph theory and make use of the tree structure of coefficient quivers of certain modules.

Recall that a graph $G$ is given by its set of vertices $V$ and a set of edges $E$ containing ordered pairs $(u,v) \in V^2$, describing an edge starting at $u$ and ending at $v$.

\begin{mydef}[{\cites{EdwardsFarr2001FragmentabilityGraphs}{EdwardsMcDiarmid1994NewUpperBoundHarmoniousColorings}}]
Let $\eps$ be a non-negative real number, and $C$ an integer. We say that a graph $G = (V, E)$ is \defined{$(C, \eps)$-fragmentable} provided there is a set $X\subseteq V$, called the fragmenting set, such that 
\begin{enumerate}
\item $\card{X}\leq \eps \card{V}$, and 
\item every component of $G\setminus X$ has at most $C$ vertices.
\end{enumerate}
Now consider a class $\Gamma$ of graphs. We will say that $\Gamma$ is \defined{$\eps$-fragmentable} provided there is an integer $C$ such that for all $G \in \Gamma$, $G$ is $(C, \eps)$-fragmentable.
Moreover, a class $\Gamma$ of graphs is called \defined{fragmentable} if \[c_f(\Gamma) \coloneqq \inf\{\eps \colon \Gamma \text{ is $\eps$-fragmentable}\} = 0.\]
\end{mydef}

\begin{remark}
We may relax the definition to say that a class $\Gamma$ of graphs is fragmentable iff for any $\eps > 0$, there are positive integers $n_0, c( \eps )$ such that if $G \in \Gamma$ is a graph with $n \geq n_0$ non-isolated vertices, then there is a set $X$ of vertices, with $\card{X} \leq \eps n$, such that each component of $G\setminus X$ has $\leq c(\eps)$ vertices.
\end{remark}

In the following, we will consider the path algebra of a given quiver $Q$. 
Recall from \cite[Section~1]{Ringel1998ExceptionalModulesTreeModules} that given a certain basis $\B$ of a representation $M$ of $Q$ (that is, a collection of basis elements from bases for vector spaces at all vertices), we say that the \emph{coefficient quiver $\Gamma(M,\B)$ of $M$ with respect to $\B$} is the quiver with vertex set $\B$ and having an arrow $b \xrightarrow{\alpha} b'$ provided the entry corresponding to $b$ and $b'$ in the matrix corresponding to $M(\alpha)$ with respect to the chosen basis $\B$ is non-zero.

We can now obtain hyperfiniteness results for a family of modules $\M$ provided a corresponding class of coefficient quivers is fragmentable.

\begin{proposition} \label{prop:TreeModulesBoundIndegreeHF}
Let $d,\ell \in \setN$.
Let $A$ be the path algebra of a quiver $Q$.
Let $\M$ be a class of indecomposable tree modules for $A$, that is, of modules $M$ such that there exists bases $\B$ of $(M_i)_{i \in Q_0}$ such that the corresponding coefficient quiver $\Gamma$ is a tree, and additionally assume that the maximal indegree of $\Gamma$ is $d$ and the maximal path length of $Q$ is $\ell$. Then $\M$ is hyperfinite.
\begin{proof}
Let $M \in \M$.
By \cite[Lemma~3.6]{EdwardsMcDiarmid1994NewUpperBoundHarmoniousColorings}, it is enough to show that the removal of at most $d^\ell$ basis elements decomposes the coefficient quiver into components of size at most half that of $M$ which are a member of $\M$. 
Since $M$ is a tree module, there is a vertex $v$ (one of the central points of the underlying tree graph) in the coefficient quiver whose removal will result in splitting the quiver into (non-connected) subtrees of size at most half that of $M$. If this vertex $v$ is a source in the coefficient quiver, it can be removed, and the induced subtrees are submodules of $M$, which are themselves tree modules in $\M$ (pick the bases given by restriction).
If $v$ is not source, at most $d$ arrows map to it. Each of their starting vertices will be removed as well, to each of which again at most $d$ arrows map. Since the path length is bounded by $\ell$, we have to remove at most $\sum_{i=1}^{\ell}d^i$ vertices to produce submodules of $M$ of dimension at most half that of $M$. These are again tree modules in $\M$.
\end{proof}
\end{proposition}

\begin{proposition}
The family of preprojective $k\Theta(d)$-modules is hyperfinite.
\begin{proof}
By Proposition~\ref{prop:AdditiveClosureStaysHyperfinite} it is enough to consider the indecomposable preprojective modules.
Now, \cite[Proposition~3]{Ringel1998ExceptionalModulesTreeModules} gives a detailed description of these modules, showing that they are tree modules and for each arrow, each basis element at the sink vertex is mapped to from at most one basis element at the source vertex. This shows that the indegree is bounded by $d$. Note that submodules of preprojective modules are preprojective.
Now apply Proposition~\ref{prop:TreeModulesBoundIndegreeHF} with degree bound $d$ and $\ell=1$ to finish the proof.
\end{proof}
\end{proposition}

In the following, we will use the sequence $a_t$ from \cite[Section~8]{Ringel1998ExceptionalModulesTreeModules} to describe the dimension vectors of the indecomposable preprojective and postinjective \mbox{$k\Theta(d)$-modules}. The sequence is defined recursively by $a_0 = 0, a_1 = 1$ and \mbox{$a_{t+1} = d a_t - a_{t-1}$} for $t \geq 1$.

\begin{lemma} \label{lemma:ClosedFormKroneckerSequence}
Fix $d \geq 3$.
Then the closed-form solution of the recurrence relation for $a_t$ is given by 
\[ a_t = \frac{\varphi^t - \psi^t}{\sqrt{d^2-4}}\, \text{ where } \varphi = \frac{d+\sqrt{d^2-4}}{2} \text{ and } \psi = \frac{d-\sqrt{d^2-4}}{2}. \]
Moreover, the quotient $a_t / a_{t+1}$ of consecutive terms converges to $\varphi^{-1}$.
\begin{proof}
Routine exercise.
\end{proof}
\end{lemma}

\begin{lemma} \label{lemma:KroneckerPIsDegreeBounds}
Fix $d \geq 2$.
Let $Q[t]$ be an indecomposable postinjective module as described in \cite[Section~8]{Ringel1998ExceptionalModulesTreeModules}, with coefficient quiver $\Gamma$ given there.
Then the outdegree of the vertices of $\Gamma$ is bounded by two, and the indegree is bounded by $(t-1)(d-2)+d$.
\begin{proof}
By the description of the arrow maps for the postinjective indecomposable module $Q[t]$ in the dual of \cite[Proposition~3]{Ringel1998ExceptionalModulesTreeModules}, the matrices of the arrows $1,\dots,d-1$ have no common non-zero columns, so the outdegree of each source with respect to the arrows $\alpha_i$, $1 \leq i \leq d-1$ is at most one. On the other hand, each row of one of these  arrow matrices contains exactly a single one.
Moreover, as the matrix for the last arrow $\alpha_d$ is constructed by concatenating zero matrices or column block matrices containing a single identity matrix block, at most one arrow $\alpha_d$ starts at each source. Indeed, the concatenation involves $t-1$ matrices---the $C(a_{j-1},a_{j})$---containing $d-2$ identity matrices, the  $E(a_{j-1})$, of varying size $a_{j-1}$ each, and one additional identity matrix $E(a_t)$.
Combining this information yields the desired result.
\end{proof}
\end{lemma}

\begin{lemma} \label{lemma:tBoundedRootDimension}
Fix $d \geq 3$. Let $M$ be a module of dimension vector $(a_{t+1},a_{t})$. Then we can express
\[t = \log_{\varphi}\left(\dim M + \sqrt{\dim M^2 + \frac{4}{d-2}}\right) - \log_{\varphi}\left(2\frac{1+\varphi}{\sqrt{d^2-4}}\right).\]
Moreover, for $\dim M \geq 3$, it holds that $t \leq c_d \sqrt{\dim M}$ for some constant $c_d$.
\begin{proof}
Clearly, $\dim M = a_{t+1} + a_{t}$. Now using the closed form of Lemma~\ref{lemma:ClosedFormKroneckerSequence}, we have that
\begin{align*}\dim M &= \frac{\varphi^{t+1} - \psi^{t+1}}{\sqrt{d^2-4}} + \frac{\varphi^t - \psi^t}{\sqrt{d^2-4}} = \frac{\varphi^{2t+1}-(\varphi\psi)^{t}\psi}{\varphi^t\sqrt{d^2-4}} + \frac{\varphi^{2t}-(\varphi \psi)^t}{\varphi^t\sqrt{d^2-4}} \\
&= \frac{\varphi^{2t+1}-\psi+\varphi^{2t}-1}{\varphi^t \sqrt{d^2-4}} = \varphi^t \frac{1+\varphi}{\sqrt{d^2-4}} - \varphi^{-t}\frac{1+\psi}{\sqrt{d^2-4}}.\end{align*}

By substitution, noting that real powers of positive numbers are positive and using $(1+\varphi)(1+\psi) = d+2$, we get
\begin{align*}t 
&= \log_{\varphi} \left( \dim M + \sqrt{\dim M^2+\frac{4}{d-2}}\right) - \log_{\varphi} \left(2\frac{1+\varphi}{\sqrt{d^2-4}}\right).\end{align*}
Now it remains to show the estimate. We first note that for $\varphi > 1$ and \[\frac{2+2\varphi}{\sqrt{d^2-4}} > \frac{d+\sqrt{d^2-4}}{\sqrt{d^2-4}} > 1,\] the subtrahend is always positive, resulting in its omittance leaving an upper bound.
Now, when $\dim M \geq 3$, 
we have \[t \leq \log_{\varphi}(1+\sqrt{2}) + \log_{\varphi}(\dim M) < \log_{\varphi}(3) + \log_{\varphi}(\dim M) \leq 2\log_{\varphi}(\dim M).\]
Now, as $\varphi > 1$, it is enough to further consider $\ln(\dim M)$.
Clearly, \[\exp\left(2\sqrt{\dim M}\right) > \frac{\left(2 \sqrt{\dim M}\right)^2}{2!} = 2\dim M > \dim M,\] so $2\sqrt{\dim M} > \ln \dim M$.
All in all, this combines to the desired inequality \[t \leq 2 \log_{\varphi}\dim M = 2 \frac{1}{\ln \varphi} \ln \dim M < \frac{4}{\ln \varphi} \sqrt{\dim M}.\]
\end{proof}
\end{lemma}

\begin{proposition}
The family of indecomposable postinjective $k\Theta(d)$-modules is hyperfinite.
\begin{proof}
We want to give a proof similar to that of \cite[Lemma~3.6]{EdwardsMcDiarmid1994NewUpperBoundHarmoniousColorings}, but adapt it to coefficient quivers of modules instead of graphs. In a first step towards proving hyperfiniteness, we hence want to find $A>0$, $0\leq \lambda < 1$ and $0<\alpha<1$ such that for any indecomposable postinjective module of dimension $n$, there are at most $An^\lambda$ basis elements that can be removed from the coefficient quiver to leave a submodule for which every indecomposable summand has dimension at most $\alpha n$, and to each summand, a similar construction can be applied, and so forth.	

Let $\eps > 0$. We put $\alpha = \frac{1}{2} + \delta$ for some $0 < \delta < \frac{1}{4}$.
Let $Q[t]$ be the indecomposable postinjective module of dimension vector $(a_{t+1},a_{t})$. Let $n = \dim Q[t]$ and assume $n>5$. We show how to split this module into small components by a sequence of stages. Before each stage $i$, all components are isomorphic to indecomposable postinjective modules, having at most $\alpha^{i-1} n$ vertices in their coefficient quivers, while the number of components with more than $\alpha^i n$ vertices is at most $\alpha^{-1}$.
Since the coefficient quiver $\Gamma$ of $Q[s]$ is a tree, there is a vertex whose removal creates subtrees of size at most $\alpha \dim Q[s]$. Note that we can assume that this vertex to remove is a sink: if the vertex to remove was a source---since all sources have outdegree at most two by Lemma~\ref{lemma:KroneckerPIsDegreeBounds} and all their neighbours are sinks---we can just remove a neighbouring sink. Note that the size of $\alpha$ allows for this modification, as we do not require that the subtrees are at most half the size of $Q[s]$.
But a removal of a sink corresponds to passing to the cokernel of an inclusion of $S_1 \embeds Q[s]$. Yet, this cokernel must have smaller dimension than $Q[s]$, and since $Q[s]$ is postinjective, must also be postinjective. This implies that the cokernel is the direct sum of indecomposable postinjective modules for smaller~$s$, as the dimension of the indecomposable postinjectives strictly increases for growing $s$. This proves that after stage $i$, all components are indecomposable postinjective modules with no more than $\alpha^i n$ vertices in their coefficient quivers.
The number of stages is the least~$k$ such that $\alpha^k n \leq L$, for an $L$ to be determined later.
Hence $\alpha^{k-1} n > L$, so that $\alpha^{1-k} < \frac{n}{L}$.

Unfortunately, this process does not create a submodule of $Q[t]$, but a factor module given by the direct sum of many smaller indecomposable postinjective modules. To attain a submodule, we must delete further vertices. Parallel to the above sequence of stages, we conduct a downstream stage to construct a submodule. In each of these stages, we only deal with the components $M$ with \[\alpha^{i} n < \dim M \leq \alpha^{i-1} n.\]
Since $\sfrac{a_t}{2} > a_{t-1}$, not in every stage a reduction takes place. But when a reduction takes place, we create submodules from $Q[s]$ by additionally removing all the vertices adjacent to the deleted sink. By the structure of the canonical coefficient quiver of $Q[s]$, there are at most  $s(d-2)+2$ such vertices, and
$s \leq c \sqrt{\dim{Q[s]}}$ by Lemma~\ref{lemma:tBoundedRootDimension}.
Note that while the dimensions of the submodules left before downstream stage~$i$ are smaller than $\dim Q[s]$, as we have removed at least one more source, the operand in this stage~$i$ is still $\alpha^{i-1} n \geq \dim Q[s]$, as we base our considerations on the original indecomposable postinjective module.
This implies that in downstream stage $i$, we remove at most $A \sqrt{\alpha^{i-1} n}$ vertices, letting $A=2+c(d-2)$. Thus, choose $\lambda = \frac{1}{2}$.
Note that in order to apply Lemma~\ref{lemma:tBoundedRootDimension} throughout, we require $\dim Q[s]\geq 3$ in all downstream stages, leading to $\alpha^{k-1} n > 2$

Now, the total number $r_i$ of vertices removed in downstream stage~$i$ is at most	
\[\alpha^{-i} A(\alpha^{i-1}n)^\lambda = A n^\lambda \alpha^{\lambda i - i - \lambda} = \frac{A n^\lambda}{\alpha} \alpha^{(1-i)(1-\lambda)}\]
and since $\alpha^{1-k} < n/L$, we have $\alpha^{1-i} < (n/L) \alpha^{k-i}$. Hence,
\begin{align*} r_i &< \frac{An^\lambda}{\alpha} (n/L)^{1-\lambda}\alpha^{(k-i)(1-\lambda)}
 = \frac{An}{\alpha L^{1-\lambda}}\beta^{k-i},\end{align*}
where $\beta = \alpha^{1-\lambda}$, with $0 < \beta < 1$. Then the total number $R$ of vertices removed from the coefficient quiver of $Q[t]$ is
\begin{align*} R = \sum_{i=1}^{k} r_i &< \frac{An}{\alpha L^{1-\lambda}} \sum_{i=1}^{k} \beta^{k-i} \\
&< \frac{An}{\alpha L^{1-\lambda}} \sum_{i \geq 0} \beta^i\\
&=\frac{An}{\alpha L^{1-\lambda}} \frac{1}{1-\beta}. \end{align*}
Since we have $1-\lambda > 0$, it follows that we can choose $L = L_\eps$ independently of~$n$ such that $R \leq \eps n$. This then shows the hyperfiniteness of $\left\{ Q[t] \colon t>0\right\}$ and thus of the postinjective component.
\end{proof}
\end{proposition}

\begin{remark}
Recall that the logarithm can be bounded above by any radical power: We have \mbox{$\ln x \leq n \sqrt[n]{x}$}.
Now the proof of the previous proposition suggests that we have an adaptation of Proposition~\ref{prop:TreeModulesBoundIndegreeHF} in the case of coefficient quivers that are graphs of genus at most $\gamma$ for fixed $\gamma \geq 0$ or for rectangular lattices of dimension $d$ for a fixed integer $d$, provided the indegree has a logarithmic bound with respect to the dimension, as these classes of graphs were shown to be fragmentable using suitable $A$, $\lambda$ and $\alpha$ (see \cite[Corollary~3.7]{EdwardsMcDiarmid1994NewUpperBoundHarmoniousColorings}).
\end{remark}

\section{A family of non-hyperfinite modules}

In the previous section we have seen that both the preprojective and the postinjective component of generalised Kronecker algebras are hyperfinite. 
Yet, Elek \cite[Section~8]{Elek2017InfiniteDimensionalRepresentationsAmenabilty} has given an argument to show that any wild Kronecker algebra is not of amenable representation type by showing that there are non-hyperfinite families of modules over the free algebras $k\langle x_1,\dots, x_r \rangle$ with $r\geq 2$ generators. Thus, the regular component must contain a non-hyperfinite family. 
We are interested in understanding and providing such a concrete counterexample of a non-hyperfinite family of modules for algebras of non-amenable representation type.

\medskip
Motivated by a similar notion of graph expanders, Barak, Impagliazzo, Shpilka and Wigderson have introduced the notion of dimension expanders (see \cites{LubotzkyZelmanov2008DimensionExpanders}{DvirShpilka2011TowardsDimensionExpandersFiniteFields}{Bourgain2009ExpandersDimensionalExpansion}{DvirWigderson2010MonotoneExpandersConstructionsApplications}).
We state a generalised notion here.
\begin{mydef}
Let $k$ be a field, $d \in \setN$, $0 < \eta \leq 1$ and $\alpha > 0$. Given a vector space $V$ and a set $\{T_1, \dots, T_d\}$ of endomorphisms of $V$, the pair $(V,\{T_i\}_{i=1}^{d})$ is called an \defined{$(\eta,\alpha)$-dimension expander of degree $d$} provided for every subspace $W \subset V$ of dimension at most $\eta \dim_k V$, we have that \[\dim_k \sum_{i=1}^d T_i(W) \geq (1+\alpha) \dim_k W.\]
\end{mydef}
\begin{remark}
As a short for $(\sfrac{1}{2},\alpha)$-dimension expanders, we may just speak of \emph{$\alpha$-dimension expanders}.
\end{remark}

Now, a sequence of dimension expanders of degree $d$ of unbounded dimensions gives rise to a non-hyperfinite family for the $d$-Kronecker algebra $k\Theta(d)$:
\begin{prop} \label{prop:DimExpGiveNonHFWildKroneckerFamily}
Let $k$ be a field, $d \in \setN$ and $\eta, \alpha > 0$. If $\{(V_i,\{T_l^{(i)}\}_{l=1}^{d})\}_{i \in I}$ is a sequence of $(\eta,\alpha)$-dimension expanders of degree $d$ such that $\dim V_i$ is unbounded, then the induced sequence of $k\Theta(d)$-modules 
\[V_i \xmanyrightarrows{T_1^{(i)}}{T_d^{(i)}} V_i\]
 is not hyperfinite.
\begin{proof}
Let $\alpha > 0$ and $\{(V_i,(T_1^{(i)},\dots,T_d^{(i)})\}$ be a sequence of $\alpha$-dimension expanders degree $d$ and of unbounded dimension $\dim V_i$.
Consider the sequence \[\left\{M_i = \left((V_i,V_i),(T_1^{(i)},\dots,T_d^{(i)})\right)\right\}_{i \in I} \in \mods k\Theta(d).\]
If this sequence was hyperfinite, for each $\eps > 0$, there exists an $L_\eps > 0$ and we can find some $M \in \{M_i \colon i \in I\}$---given by an $(\eta,\alpha)$-dimension expander space $V$---such that $\dim M = 2\dim V > 2 \frac{L_\eps}{\eta}$ with a suitable submodule $P$ exhibiting hyperfiniteness. We will denote the vector space of $P_j$ at vertex $v \in Q_0$ by $P_j(v)$. 
We have that \[\dim P_j(1) + \dim P_j(2) = \dim P_j \leq L_\eps < \eta \dim V,\] also noting that each $P_j(v)$ is a subspace of the vector space $V$ of an $(\eta,\alpha)$-dimension expander.
As each $P_j$ is a $k\Theta(d)$-module, thus $T_1(P_j(1)) + \dots + T_d(P_j(2)) \subseteq P_j(2)$, this implies that
\begin{equation} \label{eq:SinkSpaceGEQ1PlusAlphaSourceSpace} \dim P_j(2) \geq (1+\alpha) \dim P_j(1).\end{equation}
Moreover, \begin{align*}
& 2(1-\eps) \dim V  &\leq \sum_{j=1}^{t} \left(\dim P_j(1) + \dim P_j(2)\right) &\leq \sum_{j=1}^{t} \dim P_j(1) + \dim V\\
\Leftrightarrow & (1-2\eps) \dim V &&\leq \sum_{j=1}^{t} \dim P_j(1),\end{align*}
which in light of inequality \eqref{eq:SinkSpaceGEQ1PlusAlphaSourceSpace} yields that
\begin{align*} \label{eq:OneMinusEpsLEQOnePlusAlpha}
& (1-2\eps) \dim V &\leq \sum_{j=1}^{t} \frac{\dim P_j(2)}{1+\alpha} &\leq \frac{\dim V}{1+\alpha}\\
\Leftrightarrow & \phantom{(1-2} \eps &&\geq \frac{\alpha}{2(1+\alpha)},
\end{align*}
contradicting the hyperfiniteness of the sequence $\{M_i \colon i \in I\}$.
\end{proof}
\end{prop}

\begin{remark}
If all $T_i$ are such that $T_i \circ T_j = 0$ for any combination, then in general $\left(V,\{T_i\}_{i=1}^{d}\right)$ is not a dimension expander: In this situation, we have $\im T_j \subset \bigcap_{i=1}^{d} \ker T_i$ for all $1 \leq j \leq d$. Without loss of generality, we consider $0 \neq v \in \im T_1$ (if all $T_i$s were zero, the claim is obviously true). Let $W = \spann{v}$. Then $\sum_{i=1}^{d} T_i(W) = 0$, so the dimension property cannot hold for some non-trivial subspace of dimension one. Thus---unless $\eta \dim V < 1$---the pair $\left(V,\{T_i\}_{i=1}^{d}\right)$ cannot be a dimension expander.
\end{remark}

Proposition~\ref{prop:DimExpGiveNonHFWildKroneckerFamily} reduces the problem of exhibiting a non-hyperfinite family to finding families of dimension expanders for fixed $d$ and $\alpha$ such that the dimension of the vector spaces is unbounded. 
This latter question has already been asked by A.~Wigderson in 2004.
We will make use of results by \citeauthor{LubotzkyZelmanov2008DimensionExpanders} in a proposition and theorem in \cite{LubotzkyZelmanov2008DimensionExpanders} to answer it.
They provide several ways of constructing $\alpha$-dimension expanders of degree two over the complex numbers and generalise to every field of characteristic zero.

\begin{mydef}
Consider a group $\Gamma$ generated by a finite set $S$. Given a Hilbert space $H$ and a unitary representation $\rho \colon \Gamma \to U(H)$, where $U(H)$ denotes the unitary endomorphisms of $H$, the \defined{Kazhdan constant} is defined as
\[K_\Gamma^S(H,\rho) \coloneqq \inf_{0\neq v \in H} \max_{s \in S} \left\{\frac{\norm{\rho(s)v-v}}{\norm{v}} \right\}.\] 
Further, the group $\Gamma$ has \defined{property $(T)$} if
\[ K^S_\Gamma = \inf_{(H,\rho)\in \cR_0(\Gamma)} \{ K^S_\Gamma(H,\rho) \} > 0,\]
where $\cR_0(\Gamma)$ is the family of all unitary representations of $\Gamma$ which have no non-trivial $\Gamma$-fixed vector. In this case, $K^S_\Gamma$ is called the \defined{Kazhdan constant of $\Gamma$ with respect to $S$}.
\end{mydef}

This Kazhdan constant is now relevant in the following Proposition determining the expansion rate $\alpha$. In the following, by $U_n(\setC)$, we denote the group of $n\times n$ unitary matrices over $\setC$.

\begin{prop}[{\cite[Proposition~2.1]{LubotzkyZelmanov2008DimensionExpanders}}] \label{prop:LZ08DimensionExpandersFromKazhdanConstant}
Let $\rho \colon \Gamma \to U_n(\setC)$ be an irreducible unitary representation of a group $\Gamma$ with finite generating set $S$, 
then $(\setC^n,\rho(S))$ is an $\alpha$-dimension expander of degree $\card{S}$ 
where $\alpha = \frac{\kappa^2}{12}$, $\kappa = K^S_\Gamma(S\ell_n(\setC), \adj \rho)$, 
where $S\ell_n(\setC)$ denotes the subspace of all linear transformations of zero trace, 
and $\adj \rho$ is the adjoint representation on $\End(\setC^n)$ induced by conjugation.
\end{prop}

\begin{remarks}
\begin{enumerate}
\item The endomorphism space $\End(\setC^n) \isom M_n(\setC)$ and its subspace $S\ell_n(\setC)$ become Hilbert spaces via $\Euler{S,T} = \tr(ST^*)$.
\item The induced representation $\adj \rho$ on $M_n(\setC)$, given as \[\gamma \mapsto \left(T \mapsto \rho(\gamma)T\rho(\gamma)^{-1}\right),\] is unitary, as $\adj(\rho)(\gamma)$ is surjective and preserves the inner product for each $\gamma \in \Gamma$.
\item The subspace $S\ell_n(\setC)$ of trace zero matrices is invariant under $\adj \rho$, since conjugation by invertible matrices preserves the trace. Thus, $(S\ell_n(\setC),\adj \rho)$ is a unitary representation.
\item Note that if $\rho$ is irreducible, then by Schur's Lemma, $S\ell_n(\setC)$ does not have any non-trivial $\adj \rho(\Gamma)$-fixed vector:
If $T\in S\ell_n(\setC)$ is fixed by $\adj(\rho)$, then $\ker T$ is an invariant subspace of $\rho$, as $T(v) = 0$ along with $\rho(\gamma)T = T\rho(\gamma)$ implies that $T(\rho(\gamma)(v)) = 0$. 
By the irreducibility of $\rho$, $\ker T$ must be a trivial subspace. If $\ker T = 0$, we have that $T$ is invertible, even $T = \lambda \id$ for some eigenvalue $\lambda$ of $T$. But $0 = \tr T = n\lambda$, a contradiction. Thus, $\ker T = \setC^n$, so $T$ is trivial.
\end{enumerate}
\end{remarks}

In the following, we will make explicit an example using representations of $SL(2,\Z)$. This allows us to describe the Kronecker representations more easily. 
To this end, we first consider representations of the special linear group $SL(2,p)$ of $2\times 2$-matrices over the finite field of characteristic p, $\F_p$.
We recall the following two classical results, fixing some notation.
\begin{lemma} \label{lemma:IrrCplRepDimP}
For each prime $p \in \setP$, there is an irreducible, complex $p$-dimensional representation of $\SL(2,p)$.
\begin{proof}
Let $p$ be a prime number. Then $\Gamma = \SL(2,p)$ acts on $\setP_1(\F_p) = \{0, 1, \dots, p-1, \infty\}$ by
\[\pi \colon \begin{pmatrix}a & b \\ c & d\end{pmatrix} \mapsto \left( z \mapsto \frac{az+b}{cz+d}\right),\]
with the usual conventions that $\frac{x}{0} = \infty$ for $x\neq 0$ and $\frac{a\infty+b}{c\infty+d} = \frac{a}{c}$.
This permutation action extends to a permutation representation 
$\rho \colon \Gamma \to \GL_{p+1}(\setC)$,
\[g \mapsto \left(\sum_{z \in \setP_1(\F_p)} \lambda_z e_z \mapsto \sum_{z \in \setP_1(\F_p)} \lambda_z e_{\pi(z)}\right),\]
identifying $\setC^{p+1}$ with the free complex vector space on $\setP_1(\F_p)$ via \[e_1 \leftrightarrow f_0,\, \dots,\, e_p \leftrightarrow f_{p-1} \text{ and } e_{p+1} \leftrightarrow f_\infty,\]
where $\{e_1, \dots, e_{p+1}\}$ is the standard basis of $\setC^{p+1}$ and by $f_0, \dots, f_{p-1}, f_{\infty}$ we denote a standard basis of the free vector space. 
The character values of $\chi_\rho$ can be calculated via the number of fixed points of $\pi$ on representatives of the conjugacy classes of $SL(2,p)$. 
Consider the subspace $W = \{v \in \setC^{p+1} \colon \sum_{i=1}^{p+1} v_i = 0\}$ of dimension $p$. It is $\rho$-invariant and the restriction of $\rho$ to $W$ is the complement of the trivial representation in $\rho$. Using character theory, this is sufficient information to show that $\rho_{|W}$ is an irreducible complex representation of $\SL(2,p)$ (see also \cite[Section~5.2]{FultonHarris1991RepresentationTheory}).
\end{proof}
\end{lemma}

\begin{corollary} \label{cor:SL2ZIrrUnitRepOfUnbdDim}
The group $SL_2(\Z)$ has irreducible, unitary representations of unbounded dimension.
\begin{proof}
Consider the natural maps $\pi_p \colon SL_2(\Z) \to SL_2(\Z/p\Z) = SL(2,p)$ mapping each matrix to the matrix of the residue classes of its entries modulo $p$.
Let $\rho \colon SL(2,p) \to GL(V)$ be an irreducible $p$-dimensional representation. As $SL(2,p)$ is a discrete group, we can endow $V$ with an inner product in such a way to assume that $\rho$ is unitary. Now consider $\rho \circ \pi_p$. This is certainly a group homomorphism. Moreover, a subspace $W \subseteq V$ is $SL(2,p)$-invariant if and only if $W$ is $SL_2(\Z)$-invariant, showing that $\rho \circ \pi_p$ is irreducible since $\rho$ is.
\end{proof}
\end{corollary}

\begin{remark}
We may refer to the subgroups $\Gamma(p) \coloneqq \ker \left(SL_2(\Z) \to SL(2,p)\right)$ as the \emph{principal congruence subgroups} 
and have \[\Gamma(p) = \left\{\begin{pmatrix}a & b \\ c & d\end{pmatrix} \in SL_2(\Z) \colon \begin{aligned}a \equiv d \equiv 1 &\mod p\\b \equiv c \equiv 0& \mod p\end{aligned}\right\}.\]
Since the projections are surjective, the subgroups have finite index $p^3-p$ in $SL_2(\Z)$.
\end{remark}

\begin{mydef}[{\cite[Definition~4.3.1]{Lubotzky1994DiscreteGroupsExpandingGraphsInvariantMeasures}}]
Let $\Gamma$ be a finitely generated group generated by a finite symmetric set of generators $S$. Given a family $\{N_i\}_{i \in I}$ of finite index normal subgroups, $\Gamma$ is said to have \defined{property $(\tau)$ with respect to the family $\{N_i\}_{i \in I}$} provided there exists a $\kappa > 0$ such that if $(H, \rho)$ is a non-trivial unitary irreducible representation of $\Gamma$ whose kernel contains $N_i$ for some $i\in I$, then $K_{\Gamma}^S(H,\rho) > \kappa$. 
\end{mydef}

\begin{remark}
This definition is equivalent to requiring that the trivial representation is isolated in the set of all unitary representations of $\Gamma$ whose kernel contains some $N_i$ or to requiring that the non-trivial irreducible representations of $\Gamma$ factoring through $\Gamma/N_i$ for some $i\in I$ are bounded away from the trivial representation.
Further note that a finitely generated group having property $(T)$ has property $(\tau)$ with respect to all finite index normal subgroups.
\end{remark}

\begin{theorem}[{\cites[Section~1]{LubotzkyZimmer1989VariantKazhdansProperySubgroupsSemisimpleGroups}}] \label{thm:Lu94SL2ZPropertyTau}
The group $SL_2(\Z)$ has property $(\tau)$ with respect to $\{\Gamma(p)\}_{p \in \setP}$.
\begin{proof}
By Selberg's $\frac{3}{16}$ Theorem, given a congruence subgroup $\Gamma(p)$ of $SL_2(\Z)$, the smallest positive eigenvalue $\lambda_1(\Gamma(p)\backslash \setH)$ of the Laplacian on the principal modular curve $\Gamma(p)\backslash \setH$ is at least $\frac{3}{16}$. Here, $\setH$ denotes the hyperbolic plane endowed with the structure of a Riemannian manifold as in the Poincar\'{e} half-plane model. Yet, by \cite[Theorem~4.3.2]{Lubotzky1994DiscreteGroupsExpandingGraphsInvariantMeasures}, having $\lambda_1$ bound away from zero is equivalent to $SL_2(\Z)$ having property $(\tau)$ with respect to $\{\Gamma(p)\}_{p \in \setP}$.
\end{proof}
\end{theorem}

\begin{remark}
For more details and a background on the geometry, see \cite[Chapter~4]{Lubotzky1994DiscreteGroupsExpandingGraphsInvariantMeasures} or \cite[Section~3.3]{Tao2015ExpansionFiniteSimpleGroupsLieType} respectively.
\end{remark}

\begin{theorem} \label{thm:WildKroneckerNotAmenableForChar0}
Let $k$ be a field of characteristic zero. Then the wild Kronecker algebra $k\Theta(3)$ is not of amenable representation type.
\begin{proof}
By Proposition~\ref{prop:DimExpGiveNonHFWildKroneckerFamily}, it is sufficient to find a sequence of $\alpha$-dimension expanders of degree two and of unbounded dimension for some $\alpha > 0$. Now, by an application of Proposition~\ref{prop:LZ08DimensionExpandersFromKazhdanConstant}, it suffices to exhibit a sequence of irreducible, unitary representations $\rho \colon \Gamma \to U_n(\setC)$ of unbounded dimension for some group $\Gamma$ with generating set $S$ of cardinality two, such that the Kazhdan constants $K_\Gamma^S(S\ell_n(\setC),\adj \rho)$ are uniformly bounded from below by a constant $\kappa > 0$.

We let $\Gamma = \SL_2(\Z)$ with generating set $S = \{\begin{psmallmatrix}1 & 1 \\ 0 & 1\end{psmallmatrix},\begin{psmallmatrix}0 & 1 \\ -1 & 0\end{psmallmatrix}\}$. For now, we specialise to $k = \setC$.
By Corollary~\ref{cor:SL2ZIrrUnitRepOfUnbdDim}, there is a sequence $\rho_p \colon \Gamma \to U_p(\setC)$ of non-trivial irreducible, unitary representations of unbounded dimension.
Moreover, by Theorem~\ref{thm:Lu94SL2ZPropertyTau}, $\SL_2(\Z)$ has property $(\tau)$ with respect to $\{\Gamma(p)\}$, that is, there is a constant $\kappa > 0$ such that if $(H,\sigma)$ is a non-trivial unitary irreducible representation of $\SL_2(\Z)$ whose kernel contains $\Gamma(p)$ for some $p\in\setP$, then the Kazhdan constant $K_\Gamma^S(H,\sigma) > \kappa$.
Yet, by the remarks following Proposition~\ref{prop:LZ08DimensionExpandersFromKazhdanConstant}, the $(\Sl_p(\setC),\adj \rho_p)$ are unitary representations factoring through $\SL(2,p)$, that is, their kernels contain $\Gamma(p)$, and they do not contain non-trivial fixed vectors, so are irreducible.
Thus, for their Kazhdan constants we have that $K_\Gamma^S(\Sl_p(\setC),\adj \rho_p) > \kappa$.

The case for general $k$ follows as in \cite[comments after Example~3.4]{LubotzkyZelmanov2008DimensionExpanders}: Since \mbox{$\Char k = 0$}, $k$ contains~$\Q$, and the representations of Corollary~\ref{cor:SL2ZIrrUnitRepOfUnbdDim} are all defined over~$\Q$, say \mbox{$\rho_p \colon \Gamma \to \GL_p(\Q)$}. If $\card{k} \leq \aleph$, then $k$ can be embedded into $\setC$ and so can \mbox{$\GL_p(k) \subset \GL_p(\setC)$}. As the $\rho_p$ factor through a finite group, they can be unitarised over~$\setC$. We have $\setC^p = \setC \otimes_k k^p$, thus every $k$-subspace $W \subseteq k^p$ spans a $\setC$-subspace $\overline{W} \subset \setC^p$ of the same dimension. If $\rho(s) \in \GL_p(k)$, then \[\dim_k (W + \rho(s)W) = \dim_\setC(\overline{W}+\rho(s)\overline{W}).\] Since $(\setC^p,\rho_p(S))$ is a dimension expander by the above, so is $(k^p,\rho_p(S))$. 
Now, if $k$ has large cardinality and $W \subset k^n$ does not have the dimension expansion property, then the entries of a basis of $W$ generate a finitely generated field $k_1$ of characteristic zero, and we get a counterexample $W \subset k_1^n$. But $k_1^n$ is a dimension expander by the previous argument. 
\end{proof}
\end{theorem}

\begin{remark}
This proof does not use the fact that the group $\SL_2(\Z)$ has property $(\tau)$ with respect to all principal congruence subgroups, let alone all congruence subgroups. Our result follows from property $(\tau)$ with respect to infinitely many $\Gamma(p)$ such that $p$ is unbounded.
Thus, weaker versions of Selberg's Theorem should suffice in proving this. For these, see e.g., \cite[Section~3.3]{Tao2015ExpansionFiniteSimpleGroupsLieType}.
Also compare \cite[Theorem~4.4.4]{DavidoffSarnakValette2003ElementaryNumberTheoryGroupTheoryRamanujanGraphs}, where by the use of only elementary methods it is shown that the corresponding construction for graphs gives expander graphs.
\end{remark}

\begin{remark}
Put $s =\begin{psmallmatrix}1 & 1 \\ 0 & 1\end{psmallmatrix}$ and $t=\begin{psmallmatrix}0 & 1\\-1 & 0 \end{psmallmatrix}$.
Then the desired (counter)example for $k\Theta(3)$ is given by the family $\left\{ \bigl((k^p,k^p),(\id,\rho_p(s),\rho_p(t))\bigr)\right\}_{p \in \setP}$.

Considering the basis $\{e_2-e_1, \dots, e_{p}-e_{p-1}, e_{p+1}-e_{p}\}$ of $W \isom k^p$, 
we do have \[\rho_p(s) = \begin{pmatrix}0 & \dots & 0 & -1 & 1 \\ 1 & & & -1 & 1 \\ & \ddots & & \vdots & \vdots \\ & & 1 & -1 & 1 \\ 0 & \dots & 0 & 0 & 1\end{pmatrix} \in \GL_p(\Q), \quad 
\rho_3(t) = \begin{pmatrix}0 & 0 & -1\\0 &-1 & 0\\-1 & 0 & 0\end{pmatrix},\]
\[\rho_5(t) = \begin{pmatrix} 0 & 0 & 0 & 0 & -1 \\ 0 & 0 & 0 & -1 & 0 \\ 0 & -1 & 1 & -1 & 0 \\ 0 & -1 & 0 & 0 & 0 \\ -1 & 0 & 0 & 0 & 0 \end{pmatrix},\quad
\rho_7(t) = \begin{pmatrix} 0 & 0 & 0 & 0 & 0 & 0 & -1 \\ 0 & 0 & 0 & 0 & 0 & -1 & 0 \\ 0 & 0 & -1 & 1 & 0 & -1 & 0 \\ 0 & -1 & 0 & 1 & 0 & -1 & 0 \\ 0 & -1 & 0 & 1 & -1 & 0 & 0 \\ 0 & -1 & 0 & 0 & 0 & 0 & 0 \\ -1 & 0 & 0 & 0 & 0 & 0 & 0 \end{pmatrix},\]
\[\rho_{11}(t) = \begin{pmatrix} 0 & 0 & 0 & 0 & 0 & 0 & 0 & 0 & 0 & 0 & -1 \\ 0 & 0 & 0 & 0 & 0 & 0 & 0 & 0 & 0 & -1 & 0 \\ 0 & 0 & 0 & 0 & -1 & 1 & 0 & 0 & 0 & -1 & 0 \\ 0 & 0 & 0 & 0 & -1 & 1 & -1 & 1 & 0 & -1 & 0 \\ 0 & 0 & 0 & 0 & -1 & 1 & -1 & 0 & 1 & -1 & 0 \\ 0 & -1 & 1 & 0 & -1 & 1 & -1 & 0 & 1 & -1 & 0 \\ 0 & -1 & 1 & 0 & -1 & 1 & -1 & 0 & 0 & 0 & 0 \\ 0 & -1 & 0 & 1 & -1 & 1 & -1 & 0 & 0 & 0 & 0 \\ 0 & -1 & 0 & 0 & 0 & 1 & -1 & 0 & 0 & 0 & 0 \\ 0 & -1 & 0 & 0 & 0 & 0 & 0 & 0 & 0 & 0 & 0 \\ -1 & 0 & 0 & 0 & 0 & 0 & 0 & 0 & 0 & 0 & 0 \end{pmatrix}.\]
While the latter matrices $\rho_p(t)$ share a pattern, the strict rule to construct them ad-hoc is unclear. Note that we do have a certain symmetry $a_{i,j} = a_{p+1-i,p+1-j}$.
\end{remark}

\begin{theorem}
Let $k$ be any field. Then there exists some $d\geq3$ such that $k\Theta(d)$ is not of amenable representation type.
\begin{proof}
Similarly to the proof of Theorem~\ref{thm:WildKroneckerNotAmenableForChar0}, we must find a family of dimension expanders of unbounded dimension.
From \cite[Theorem~1.2]{DvirWigderson2010MonotoneExpandersConstructionsApplications}, we know that to give an explicit construction of degree-$d$ dimension expanders, it is sufficient to have an explicit construction of so-called $d$-monotone expander graphs. Note that the authors attribute this to implicit work in the initial publication of \cite{DvirShpilka2011TowardsDimensionExpandersFiniteFields}. An explicit construction of a constant-degree (discrete) monotone expander graph was suggested and outlined in \cite{Bourgain2009ExpandersDimensionalExpansion} and presented in \cite[Corollary~2]{BourgainYehudayoff2013ExpansionSL2RMonotoneExpanders}.

Now, given any field $k$, this construction allows us to find degree-$d$ dimension expanders of arbitrarily large dimension, thus showing that the wild $d$-Kronecker algebras $k\Theta(d)$ are not of amenable representation type.
\end{proof}
\end{theorem}

\section{Finitely controlled wild algebras}
Let us now consider a weaker version of hyperfiniteness and how it can be used to prove the non-amenability of a large class of wild algebras.

\begin{mydef}
Let $k$ be a field, let $A$ be a finite dimensional $k$-algebra and let \mbox{$\M \subseteq \mods A$} be a family of finite dimensional $A$-modules. $\M$ is called \defined{weakly hyperfinite}\index{hyperfinite!weakly} provided for every $\eps >0$ there exists some $L_\eps > 0$ such that for every \mbox{$M \in \M$} there is a homomorphism $\vartheta \colon N \to M$ for some $N \in \mods A$ such that
\begin{equation} 
\begin{split}
\dim_k \ker \vartheta \leq \eps \dim M, \\
\dim_k \coker \vartheta \leq \eps \dim M, 
\end{split}\end{equation}
and there are modules $N_1,\dots,N_t \in \mods A$ with $\dim_k N_i \leq L_\eps$ such that $N \isom \bigoplus_{i=1}^{t}N_i$.

A $k$-algebra $A$ is said to be of \defined{weakly amenable representation type}\index{representation type!amenable!weakly} provided $\mods A$ itself is a weakly hyperfinite family.
\end{mydef}

We see that the term ``weakly hyperfinite'' is suitably chosen:

\begin{prop} \label{prop:HFimplieswHF}
Let $A$ be a finite dimensional $k$-algebra. If $\M \subseteq \mods A$ is hyperfinite, then $\M$ is weakly hyperfinite.
\begin{proof}
Let $\eps > 0$. By the hyperfiniteness, there is some $L_\eps > 0$.
Let $M \in \M$. Then there is $N \subseteq M$ with $\dim N \geq (1-\eps) \dim M$ and $N = \bigoplus_{i=1}^{t} N_i$ with $\dim N_i \leq L_\eps$.
Let~$\vartheta$ be the inclusion of the submodule $N \embeds M$. 
Then $\ker \vartheta = 0$, and $\coker \vartheta \isom M / N$, thus $\dim \coker \vartheta = \dim M - \dim N \leq \eps \dim M$.
We have shown weak hyperfiniteness with the same $\eps$ and $L_\eps$.
\end{proof}
\end{prop}

We next turn to the relation between dimension expanders and examples of non-hyperfinite families and generalise this result to families which are not weakly hyperfinite.

\begin{prop} \label{prop:DimExpGiveNonwHFWildKroneckerFamily}
Let $k$ be a field, $d \in \N$ and $\eta, \alpha > 0$. If $\{(V_i,\{T_l^{(i)}\}_{l=1}^{d})\}_{i \in I}$ is a family of $(\eta,\alpha)$-dimension expanders of degree $d$ such that $\dim V_i$ is unbounded, then the induced family of $k\Theta(d)$-modules 
\[V_i \xmanyrightarrows{T_1^{(i)}}{T_d^{(i)}} V_i\]
is not weakly hyperfinite.
\begin{proof}
Assume to the contrary that the induced family $\M = \{M_i \colon i \in I\}$ was weakly hyperfinite. Let $\eps > 0$. Then there exists $L_\eps > 0$ such that for all $M \in \M$ there exists $\vartheta \colon P \to M$ such that $\dim \ker \vartheta \leq \eps \dim M$, $\dim \coker \vartheta \leq \eps \dim M$ and $P = \bigoplus_{j=1}^{s} P_j$ where $\dim P_j \leq L_\eps$.
For some $(\alpha,\eta)$-dimension expander $(V,\{T_\ell\})$ of degree $d$, such that $\dim V \geq \frac{1}{\eta} L_\eps$,
consider \[M = \bigl((V,V),(T_1,\dots,T_d)\bigr).\]
Letting $\vartheta_j \colon P_j \xrightarrow{\iota_j} \bigoplus P_\ell \xrightarrow{\vartheta} M$, we see that $\left(\vartheta_j(P_j)\right)(1)$ is a subspace of $V$ and $\dim \left(\vartheta_j(P_j)\right)(1) < \eta \dim V$.
Moreover, $\vartheta_j(P_j)$ is a $k\Theta(d)$-module. It follows 
via the expander property that
\begin{equation} \label{eq:DimExpPropOnSinkSpace}
\dim\left(\vartheta_j(P_j)\right)(2) \geq \dim \sum_{\ell=1}^{d} T_\ell\left((\vartheta_j(P_j))(1)\right) \geq (1+\alpha) \dim \left(\vartheta_j(P_j)\right)(1).\end{equation}
As $0 \to \im \vartheta \to M \to \coker \vartheta \to 0$ is exact, we have that
\[\dim \im \vartheta = \dim M - \dim \coker \vartheta \geq (1-\eps) \dim M = 2(1-\eps) \dim V.\]
On the other hand, 
\begin{align*}
\dim \im \vartheta &\leq \sum_{j=1}^{s} \bigl[\dim \left(\vartheta_j(P_j)\right)(1) + \dim \left(\vartheta_j(P_j)\right)(2)\bigr] \\
&\stackrel{\eqref{eq:DimExpPropOnSinkSpace}}{\leq} \sum_{j=1}^{s} \left[\left(1+\frac{1}{1+\alpha}\right) \dim \left(\vartheta_j(P_j)\right)(2)\right] = \left(\frac{2+\alpha}{1+\alpha}\right) \sum_{j=1}^{s} \dim \left(\vartheta_j(P_j)\right)(2).
\end{align*}
Next, note that \[0 \to \ker \vartheta(2) \embeds \left(\bigoplus_{j=1}^{s} P_j\right)(2) \xrightarrow{\begin{psmallmatrix}\vartheta_1(2) & \dots & \vartheta_s(2)\end{psmallmatrix}} V\] is exact,
showing that \[\dim \bigoplus_{j=1}^{s} P_j(2) \leq \dim \ker \vartheta(2) + \dim V.\]
Now, compare these to estimates to get
\begin{align*}
2(1-\eps) \dim V &\leq \dim \im \vartheta \leq \left(\frac{2+\alpha}{1+\alpha}\right) \sum_{j=1}^{s} \dim \vartheta_j(P_j)(2) \leq \left(\frac{2+\alpha}{1+\alpha}\right) \sum_{j=1}^{s} \dim P_j(2)\\
& \leq \left(\frac{2+\alpha}{1+\alpha}\right) \left(\dim V + \dim \ker \vartheta(2)\right) \leq \left(\frac{2+\alpha}{1+\alpha}\right) (1+2\eps) \dim V\\
\Leftrightarrow 2-2\eps \phantom{\dim V}&\leq \left(\frac{2+\alpha}{1+\alpha}\right) (1+2\eps) \Leftrightarrow 2-\left(\frac{2+\alpha}{1+\alpha}\right)(1+2\eps) \leq 2\eps \\
\Leftrightarrow \eps \phantom{\dim V}&\geq \frac{2+2\alpha-(2+\alpha)(1+2\eps)}{2(1+\alpha)} = \frac{\alpha-4\eps - 2\eps \alpha}{2(1+\alpha)}\\
\Leftrightarrow \frac{\alpha}{2(1+\alpha)} &\leq \eps + \frac{2\eps+\eps\alpha}{1+\alpha} = \eps\left(1+\frac{2+\alpha}{1+\alpha}\right) = \eps \left(\frac{1+\alpha+2+\alpha}{1+\alpha}\right)\\
\Leftrightarrow \eps \phantom{\dim V}&\geq \frac{\alpha}{2(1+\alpha)}\frac{1+\alpha}{3+2\alpha} = \frac{\alpha}{6+4\alpha} > 0,
\end{align*}
contradicting the weak hyperfiniteness of $\M$.
\end{proof}
\end{prop}

\begin{corollary} \label{cor:SomeKThetaDNotweaklyAmenable}
Let $k$ be any field. Then there exists some $d \geq 3$ such that $k\Theta(d)$ is not of weak amenable representation type.
\end{corollary}

We recall a notion of wildness originally due to \textcite{Ringel2002CombinatorialRepresentationTheoryHistoryFuture} that is suitable for our purpose.

\begin{mydef}[{\cites[Definition~2.1]{Han2001ControlledWildAlgebras}[Section~4]{GregoryPrest2016RepresentationEmbeddingsInterpretationFunctorsControlledWildAlgebras}}]
An algebra $A$ is \defined{controlled wild} and \defined{controlled by $\C$} provided there exist a faithful exact functor \[F\colon \mods k\langle x,y\rangle \to \mods A,\] and a full subcategory $\C$ of $\mods A$ which is closed under direct sums and direct summands such that for any $X$ and $Y$ in $\mods k\langle x,y\rangle$,
\[\Hom_A\bigl(F(X),F(Y)\bigr) = F\left(\Hom_{k\langle x,y\rangle} (X,Y)\right) \oplus \Hom_A(F(X),F(Y))_\C,\]
and
\[\Hom_A(F(X),F(Y))_\C \subseteq \rad \Hom_A(F(X),F(Y)).\]
What is more, we say that $A$ is \defined{finitely controlled wild}\index{wild!finitely controlled} if it is controlled by $\add(C)$ for some $C \in \mods A$.
\end{mydef}

Here, $\Hom_A(X,Y)_\C$ denotes the set of those $A$-homomorphisms  $X \to Y$ factoring through $\C$.

\begin{remark}

\cite[Problem~10]{Ringel2002CombinatorialRepresentationTheoryHistoryFuture} asks whether all wild algebras are controlled wild. This result was announced by Y.~Drozd in 2007 (even for finitely controlled), but has not yet been published. 
\end{remark}

Finally, we can conclude this note by showing that finitely controlled wild algebras are not weakly amenable. 

\begin{theorem} \label{thm:FinCtrWildNotWeaklyAmenable}
Let $k = \bar{k}$ be an algebraically closed field and let $A$ be a finite dimensional $k$-algebra. If $A$ is finitely controlled wild, then $A$ is not of weakly amenable representation type.
\begin{proof}
Let $d$ be as in Corollary~\ref{cor:SomeKThetaDNotweaklyAmenable}.
Since $A$ is finitely controlled wild, by \cite[Lemma~2.4]{Han2001ControlledWildAlgebras} there is a faithful and exact functor $F \colon \mods k\Theta(d) \to \mods A$, which is a finitely controlled representation embedding in the sense of \cite[Section~4]{GregoryPrest2016RepresentationEmbeddingsInterpretationFunctorsControlledWildAlgebras} controlled by a full subcategory $\C = \add(C)$ for some $C \in \mods A$. Then by \cite[Theorem~4.2]{GregoryPrest2016RepresentationEmbeddingsInterpretationFunctorsControlledWildAlgebras}, there is a functor $G \colon \mods A \to \mods k\Theta(d)$ such that $(G \circ F) (M) \isom M$ for all $M \in \mods k\Theta(d)$. 
Let us denote by $K = {_{k\Theta(d)} k\Theta(d)}$ the Kronecker algebra as a left module over itself.
Now, the functor $G$ is given on the objects by
\[G(X) = \faktor{\Hom_A(F(K),X)}{\Hom_A(F(K),X)_\C.}\]
We can find a $\C$-preenvelope of $F(K)$,
\[\Delta \colon F(K) \to \C_{F(K)} \isom C^n,\]
where $n = \dim_k \Hom_A \left(F(K),C\right)$, and we have
\[\Hom(\Delta,-) \colon \Hom_A(C^n,-) \to \Hom_A(F(K),-), \quad f \mapsto f \circ \Delta.\]
We note that $\Hom_A(\Delta,C')$ is surjective for each $C' \in \C$ and that every morphism in the image factors through $\C$. 
It follows that $G$ is the cokernel functor of $\Hom_A(\Delta,-)$, and for all $X \in \mods A$ we have that
\begin{equation} \label{eq:IntFunctorIsCokernel} 0 \to \Hom_A(F(K),X)_\C = \im \Hom_A(\Delta,X) \embeds \Hom_A(F(K),X) \onto G(X) \to 0.\end{equation}
Assume that $\mods A$ was weakly hyperfinite. Then, for every $\tilde{\eps} > 0$, there exists ${L}^{\mods A}_{\tilde{\eps}} > 0$ fulfilling the usual conditions. Hence, given $M \in \mods k\Theta(d)$, we can find $\vartheta \colon N \to F(M)$ in $\mods A$ such that $N \isom \bigoplus_{i=1}^{s} N_i$, with $\dim N_i \leq {L}^{\mods A}_{\tilde{\eps}}$ and $\ker \vartheta, \coker \vartheta \leq \tilde{\eps} \dim F(M)$.
We shall consider the exact sequences
\begin{align*}
\epsilon \colon 0 \to \ker \vartheta \xrightarrow{\alpha} N \xrightarrow{\beta} \im \vartheta \to 0,\\
\eta \colon 0 \to \im \vartheta \xrightarrow{\gamma} M \xrightarrow{\delta} \coker \vartheta \to 0.
\end{align*}

Connecting two sequences of type \eqref{eq:IntFunctorIsCokernel} by $b = \Hom(F(K),\vartheta)$, via an application of the Snake Lemma, we get the following commutative diagram.
\begin{center}
\begin{tikzpicture}
\matrix (m) [matrix of math nodes,nodes={anchor=center},row sep=2em,column sep=3em] 
  {
  0 & \ker a  & \ker b & \ker G(\vartheta) \\
  0 & \im (\Delta,N) & {(F(K),N)} & G(N) & 0\phantom{.} \\
  0 & \im (\Delta,F(M)) & {(F(K),F(M))} & G(F(M)) & 0\phantom{.} \\
  & \coker a & \coker b & \coker G(\vartheta) & 0.\\
     };
  \path[->]
    (m-2-1) edge (m-2-2)
	(m-2-2) edge [right hook->] node [above] {$f$} (m-2-3)
			edge node [right] {$a$} (m-3-2)
	(m-2-3) edge [->>] node [above] {$g$} (m-2-4)
			edge node [right] {$b$} (m-3-3)
	(m-2-4) edge (m-2-5)
			edge node [right] {$G(\vartheta)$} (m-3-4)
	(m-3-1) edge (m-3-2)
	(m-3-2) edge [right hook->] node [above] {$f'$} (m-3-3)
	(m-3-3) edge [->>] node [above] {$g'$} (m-3-4)
	(m-3-4) edge (m-3-5);
  \path[->,gray]
    (m-1-2) edge (m-2-2)
	(m-1-3) edge (m-2-3)
	(m-1-4) edge (m-2-4)
	(m-3-4) edge (m-4-4)
	(m-3-2) edge (m-4-2)
	(m-3-3) edge (m-4-3);
  \draw[->,dashed]
	(m-1-1) edge (m-1-2)
	(m-1-2) edge (m-1-3)
	(m-1-3) edge (m-1-4)
	(m-4-2) edge (m-4-3)
	(m-4-3) edge (m-4-4)
	(m-4-4) edge (m-4-5);
  \draw[->,dashed,rounded corners] 
	(m-1-4) -| ($(m-2-5.east)+(.5,0)$) node [right] {$d$} |- ($(m-3-3)!0.65!(m-2-3)$) -| ($(m-3-1.west)+(-.5,0)$) |- (m-4-2);
\end{tikzpicture}
\end{center}
We see that 
$\dim \ker G(\vartheta) = \dim \ker d + \dim \im d \leq \dim \ker b + \dim \coker a$ 
and \mbox{$\dim \coker G(\vartheta) \leq \dim \coker b$}. 

Given $\varphi \in \ker b$, 
there is a unique \mbox{$\varphi'\colon F(K) \to \ker \varphi$} such that \mbox{$\Ker \varphi \circ \varphi' = \varphi$}, showing that \mbox{$\Hom(F(K), \ker \varphi) \onto \ker b$}, hence $\dim \ker b \leq \dim \Hom\left(F(K),\ker \varphi\right)$.

Since $\im b \subseteq \im \Hom(F(K),\gamma)$, $\coker b$ has as submodule $\im \Hom(F(K),\gamma) / \im b$,
and the quotient considered as a vector space is
$\coker \Hom_A(F(K),\gamma)$. 
It now follows that
\begin{align*}
\dim \text{quotient} &= \dim (F(K),F(M)) - \dim \im (F(K),\gamma) \\
&= \dim (F(K),F(M)) - \dim \ker (F(K),\delta)\\
&= \dim \im (F(K),\delta) \leq \dim \Hom_A(F(K),\coker \vartheta),\end{align*}
using the left exactness of $\Hom_A(F(K),-)$.
On the other hand, as we recall that \mbox{$\epsilon^* \colon (F(K),\im \vartheta) \to {^{1}}(F(K),\ker \vartheta)$}, 
\begin{align*}
\dim \text{submodule} &= \dim \im (F(K),\gamma) - \dim \im b \\
&= \dim \im (F(K),\gamma) - \dim \im \big((F(K),\gamma) \circ (F(K),\beta)\bigr)\\
&= \dim (F(K),\im \vartheta) - \dim \im (F(K),\beta) = \dim (F(K),\im \vartheta) - \dim \ker \epsilon^* \\
&= \dim \im \epsilon^* \leq \dim \Ext{A}{1}(F(K),\ker \vartheta),\end{align*}
as $\Hom(F(K),\gamma)$ is a monomorphism. 
Combining these two inequalities shows that 
\mbox{$\dim \coker b \leq \dim \Hom_A(F(K), \coker \vartheta) + \dim \Ext{A}{1}(F(K),\ker \vartheta)$}. 

We are left to deal with $\coker a$. Here, we consider yet another diagram with exact rows which we complete again by an application of the Snake Lemma.
\begin{center}
\begin{tikzpicture}
\matrix (m) [matrix of math nodes,nodes={anchor=center},row sep=2em,column sep=3em] 
  {
  0 & \ker \sigma  & \ker \tau& \ker a \\
  0 & \ker (\Delta,N) & {(C^n,N)} & \im (\Delta,N) & 0\phantom{.} \\
  0 & \ker (\Delta,F(M)) & {(C^n,F(M))} & \im (\Delta,F(M)) & 0\phantom{.} \\
  & \coker \sigma & \coker \tau & \coker a & 0.\\
     };
  \path[->]
    (m-2-1) edge (m-2-2)
	(m-2-2) edge [right hook->] (m-2-3)
			edge node [right] {$\sigma$} (m-3-2)
	(m-2-3) edge [->>] (m-2-4)
			edge node [right] {$\tau = (C^n,\vartheta)$} (m-3-3)
	(m-2-4) edge (m-2-5)
			edge node [right] {$a$} (m-3-4)
	(m-3-1) edge (m-3-2)
	(m-3-2) edge [right hook->] (m-3-3)
	(m-3-3) edge [->>] (m-3-4)
	(m-3-4) edge (m-3-5);
\path[->,dashed]
    (m-1-1) edge (m-1-2)
    (m-1-2) edge (m-1-3)
	(m-1-3) edge (m-1-4)
	(m-4-2) edge (m-4-3)
	(m-4-3) edge (m-4-4)
	(m-4-4) edge (m-4-5);
  \draw[->,gray]
	(m-1-2) edge (m-2-2)
	(m-1-3) edge (m-2-3)
	(m-1-4) edge (m-2-4)
	(m-3-2) edge (m-4-2)
	(m-3-3) edge (m-4-3)
	(m-3-4) edge (m-4-4);
  \draw[->,dashed,rounded corners] 
	(m-1-4) -| ($(m-2-5.east)+(.5,0)$) |- ($(m-3-3)!0.65!(m-2-3)$) -| ($(m-3-1.west)+(-.5,0)$) |- (m-4-2);
\end{tikzpicture}
\end{center}
As above, $\dim \coker a \leq \dim \coker \tau$. 
By the same line of argument as before, we have \mbox{$\dim \coker \tau \leq \dim \Hom_A(C^n,\coker \vartheta) + \dim \Ext{A}{1}(C^n,\ker \vartheta)$}. 

Now consider that given $X \in \mods A$, there are $m \in \setN$ and $Y_X = \ker (A^m \onto X)$ such that $0 \to Y_X \to A^m \to X \to 0$ is exact.
Applying $\Hom_A(-,\ker \vartheta)$ 
allows us to deduce that 
$\dim \Ext{A}{1}(X,\ker \vartheta) 
\leq (\dim X + \dim Y_X) \dim \ker \vartheta$.
All in all, we have
\begin{align*}
\dim \coker G(\vartheta) &\leq \dim \coker b \leq \dim \Hom_A(F(K),\coker \vartheta) + \dim \Ext{A}{1}(F(K),\ker \vartheta)\\
&\leq \dim F(K) \tilde{\eps} \dim F(M) + \left(\dim Y_{F(K)} + \dim F(K)\right) \tilde{\eps} \dim F(M)\\
&\leq \dim F(K) \left(2\dim F(K)+\dim Y_{F(K)}\right) \tilde{\eps} \dim M.
\end{align*}
On the other hand,
\begin{align*}
\dim \ker G(\vartheta) &\leq \ker b + \dim \coker a \\
&\leq \dim \Hom_A(F(K),\ker \vartheta) + \dim \Hom_A(C^n,\coker \vartheta)\\
&\phantom{\leq \dim \Hom_A(F(K),\ker \vartheta)~} + \dim \Ext{A}{1}(C^n,\ker \vartheta)\\
&\leq \left(\dim F(K) + 2n \dim C + \dim Y_{C^{n}}\right) \tilde{\eps} \dim F(M)\\
&\leq \dim F(K) \left(\dim F(K) + 2 \dim F(K) (\dim C)^2 + \dim Y_{C^n}\right) \tilde{\eps} \dim M.
\end{align*}

\noindent To conclude the proof, we now choose \[\tilde{\eps} = \eps{\left(\dim F(K) \cdot \left(\dim F(K) + 2\dim F(K) (\dim C)^2 + \dim Y_{F(K)} + \dim Y_{C^n}\right)\right)}^{-1}\] and put $L_\eps = \dim F(K) {L}^{\mods A}_{\tilde{\eps}}$. Note that $\tilde{\eps}$ depends only on properties of $A$ and its controlled wildness. 
Then we have \[G(\vartheta) \colon G(N) \to G(F(M)) \isom M,\]
such that $G(N) \isom \bigoplus_{i=1}^{s} G(N_i)$ by the additivity of $G$, with \begin{align*} \dim_k G(N_i) &\leq \dim_k \Hom_A(F(K),N_i) - \dim_k \Hom_A(F(K),N_i)_\C \\ &\leq \dim_k F(K) \dim_k N_i \leq L_\eps,\end{align*}
and \[\dim_k \ker G(\vartheta), \, \dim_k \coker G(\vartheta) \leq \eps \dim_k M.\]
This shows that $\mods k\Theta(d)$ is weakly hyperfinite, a contradiction to Corollary~\ref{cor:SomeKThetaDNotweaklyAmenable}.
Hence, $\mods A$ cannot be weakly hyperfinite, so $A$ is not of weakly amenable representation type.
\end{proof}
\end{theorem}

\begin{corollary}
Let $k = \bar{k}$ be an algebraically closed field and let $A$ be a finite dimensional $k$-algebra.
If $A$ is a finitely controlled wild, then $A$ is not of amenable representation type.
\end{corollary}

\begin{acknowledgement}
These notes are based on work done during the author’s doctorate studies at Bielefeld University. The author would like to thank his supervisor Professor W. Crawley-Boevey for his advice and guidance.
\end{acknowledgement}

{
 \printbibliography
}
\end{document}